\documentclass[11pt]{article}
\usepackage{amsmath,amssymb}
\usepackage[latin1]{inputenc}
\usepackage[dvips]{graphicx}

\topskip=\baselineskip  \headsep=0in
\newtheorem{theorem}{Theorem}[section]
\newtheorem{lemma}[theorem]{Lemma}
\newtheorem{proposition}[theorem]{Proposition}
\newtheorem{corollary}[theorem]{Corollary}

\newtheorem{example}[theorem]{Example}

 \newcommand{\qed}{\enspace\vrule  height6pt  width4pt  depth2pt}

 \newenvironment{proof}{\par\noindent{\bf Proof.}}{$\qed$\par\bigskip}

 \newcommand{\Aut}{\mbox{\rm Aut}}

 \newcommand{\inv}{^{-1}}
 \newcommand{\N}{{\mathbb N}}
 \newcommand{\Z}{{\mathbb Z}}

 \newcommand{\C}{{\mathbb C}}

 \newcommand{\G}{{\cal G}}

 \newcommand{\matriz}[1]{\begin{array} #1 \end{array}}
 
 \newcommand{\GEN}[1]{\langle #1 \rangle}

 \newcommand{\Sym}{\mbox{\rm Sym}}
 
 \newcommand{\Fa}{\mbox{\rm Fa}}
 \newcommand{\FaM}{\mbox{\rm FaM}}
 \newcommand{\id}{\mbox{\rm id}}

\begin{document}

 \title{Involutive Yang-Baxter Groups\thanks{The authors have been partially supported by grants of
    DGI MEC-FEDER (Spain) MTM2005-00934, Generalitat de Catalunya 2005SGR00206;
    Onderzoeksraad of Vrije Universiteit Brussel,
    Fonds voor Wetenschappelijk Onderzoek (Belgium),
    Flemish-Polish bilateral agreement  BIL2005/VUB/06;
    DGI MEC-FEDER (Spain) MTM2006-06865 and Fundaci\'{o}n S\'{e}neca of Murcia 04555/GERM/06.
    \newline 2000 {\em Mathematics Subject Classification.} Primary 81R50, 20F29, 20B35, 20F16.
    \newline {\em Keywords:} Yang-Baxter equation, involutive non-degenerate solutions, group of $I$-type,
    finite solvable group.}}
 \author{Ferran Ced\'{o}, Eric Jespers and \'{A}ngel del R\'{\i}o}
 \date{}
 \maketitle

\begin{abstract}
In 1992 Drinfeld posed the question of finding the set theoretic solutions of the Yang-Baxter
equation. Recently, Gateva-Ivanova and Van den Bergh and Etingof, Schedler and Soloviev have shown
a group theoretical interpretation of involutive non-degenerate
solutions.
Namely, there is a one-to-one correspondence between
involutive non-degenerate
solutions
on finite sets and groups of $I$-type. A group $\G$ of $I$-type is a group isomorphic to a subgroup
of $\Fa_n\rtimes \Sym_n$ so that the projection onto the first component is a bijective map, where
$\Fa_n$ is the free abelian group of rank $n$ and $\Sym_{n}$ is the symmetric group of degree $n$.
The projection of $\G$ onto the second component $\Sym_n$ we call an involutive Yang-Baxter
group (IYB group). This suggests the following strategy to attack Drinfeld's problem for
involutive non-degenerate set theoretic solutions. First classify the IYB groups and second, for a
given IYB group $G$, classify the groups of $I$-type with $G$ as associated IYB group.
It is known that every IYB group is solvable. In this paper some results supporting the converse
of this property are obtained. More precisely, we show that some classes of groups are IYB groups.
We also give a non-obvious method to construct infinitely many groups of $I$-type (and hence
infinitely many involutive non-degenerate set theoretic solutions of the Yang-Baxter equation) with
a prescribed associated IYB group.
\end{abstract}

\section{Introduction}

In a paper on statistical mechanics by Yang
\cite{Yang}, the quantum Yang-Baxter equation
appeared. It turned out to be one of the basic
equations in mathematical physics and  it lies at
the foundation of the theory of quantum groups. One
of the important unsolved   problems is to discover
all the solutions $R$ of the quantum Yang-Baxter
equation (note that $R:V\otimes V \rightarrow
V\otimes V$, with $V$ a vector space). In recent
years, many solutions have been found and the
related algebraic structures have been intensively
studied (see for example \cite{kassel}). Drinfeld,
in \cite{drinfeld}, posed the question of finding
the simplest solutions, that is, the solutions $R$
that are induced by a bijective mapping $\mathcal{R}
: X\times X \rightarrow X\times X$, where $X$ is a
basis for $V$. Let $\tau : X^{2}\rightarrow X^{2}$
be the map defined by $\tau (x,y)=(y,x)$. In
\cite{ESS} it is shown that $\mathcal{R}$ yields a
solution if and only if the bijective mapping
$r=\tau \circ \mathcal{R}$ satisfies
$r_{1}r_{2}r_{1}=r_{2}r_{1}r_{2}$, where
$r_1=r\times \id_X\colon X^3\rightarrow X^3$ and
$r_2=\id_X\times r\colon X^3\rightarrow X^3$.  In
this case, one says that $r$ is a set theoretic
solution of the Yang-Baxter equation.

Set theoretic solutions $r:X^{2}\rightarrow X^{2}$
of the Yang-Baxter equation (with $X$ a finite set)
that are involutive (i.e., $r^{2}$ is the identity
map on $X^{2}$) and that are (left) non-degenerate
recently have received a lot of attention by
Etingof, Schedler and Soloviev \cite{ESS},
Gateva-Ivanova and Van den Bergh \cite{Gat,GIVdB},
Rump \cite{Rump}, Jespers and Okninski
\cite{JO,JObook} and others. Recall that a bijective
map $$\begin{array}{cccc} r\colon & X\times
X&\longrightarrow &X\times X\\
&(x,y)&\mapsto&(f_x(y), g_y(x)) \end{array}$$  is
said to be right (respectively, left) non-degenerate
if each map $f_{x}$ (respectively, $g_{x}$) is
bijective. Note that, since $X$ is finite, one can
show that an involutive set theoretic solution $r$
of the Yang-Baxter equation is right non-degenerate
if and only if it is left non-degenerate (see
\cite[Theorem~1.3]{GIVdB}, \cite[Corollary~2.3]{JO}
and \cite[Corollary 8.2.4]{JObook}).

Gateva-Ivanova and Van den Bergh in \cite{GIVdB},
and Etingof, Schedler and Soloviev in \cite{ESS},
gave a beautiful group theoretical interpretation of
involutive non-degenerate set theoretic solutions of
the Yang-Baxter equation. In order to state this, we
need to introduce some notation. Let $\FaM_n$ be the
free abelian monoid of rank $n$ with basis
$u_1,\dots , u_n$. A monoid $S$ generated by a set
$X=\{ x_1,\dots , x_n\}$ is said to be of left
$I$-type if there exists a bijection (called a left
$I$-structure) $v\colon \FaM_n\longrightarrow S$
such that $v(1)=1$ and $\{ v(u_1a),\dots
,v(u_na)\}=\{ x_1v(a),\dots ,x_nv(a)\}$, for all
$a\in \FaM_n$. In \cite{GIVdB} it is shown that
these monoids $S$ have a presentation $$S=\langle
x_1,\dots ,x_n\mid x_ix_j=x_kx_l\rangle ,$$ with
$n\choose 2$ defining relations so that every word
$x_ix_j$ with $1\leq i,j\leq n$ appears at most once
in one of the relations. Such presentation induces a
bijective map $r\colon X\times X\longrightarrow
X\times X$ defined by
    $$r(x_i,x_j)=
    \left\{\matriz{{ll}
    (x_k,x_l), & \text{if }x_ix_j=x_kx_l \text{ is a defining relation for } S;\\
    (x_i,x_j), & \text{otherwise.}
    }\right.$$
Furthermore, $r$ is an involutive left
non-degenerate set theoretical solution of the
Yang-Baxter equation. Conversely, for every
involutive left non-degenerate set theoretical
solution of the Yang-Baxter equation $r\colon
X\times X\longrightarrow X\times X$ and every
bijection $v:\{u_1,\dots,u_n\}\rightarrow X$ there
is a unique left $I$-structure $v:\FaM_n\rightarrow
S$ extending $v$, where $S$ is the semigroup given
by the following presentation $S=\GEN{X \mid ab=cd,
\text{ if } r(a,b)=(c,d)}$
(\cite[Theorem~8.1.4.]{JObook}). Furthermore, it is
proved in \cite{GIVdB} that a monoid of left
$I$-type has a group of fractions, which is called a
group of left $I$-type.

In \cite{JO} Jespers and Okni\'{n}ski proved that a
monoid $S$ is of left $I$-type if and only if it is
of right $I$-type and obtained an alternative
description of monoids and groups of $I$-type.
Namely, it is shown that a monoid is of $I$-type if
and only if it is isomorphic to a submonoid $S$ of
the semidirect product $\FaM_n\rtimes \Sym_n$, with
the natural action of $\Sym_n$ on $\FaM_n$ (that is,
$\sigma (u_i)=u_{\sigma(i)}$), so that the
projection onto the first component is a bijective
map, that is
    $$S=\{ (a,\phi (a))\mid a\in \FaM_n\},$$
for some map $\phi\colon \FaM_n\rightarrow \Sym_n$.
In that case the map $\phi$ extends uniquely to a
map $\phi\colon\Fa_n\rightarrow \Sym_n$, where
$\Fa_n$ is the free abelian group of rank $n$, and
the corresponding group of $I$-type $SS^{-1}$ is
isomorphic to a subgroup ${\cal G}$ of the
semidirect product $\Fa_n\rtimes \Sym_n$ so that the
projection onto the first component is a bijective
map, that is
    \begin{eqnarray}\label{Itypegroup}
    {\cal G}=\{ (a,\phi(a))\mid a\in \Fa_n\}.
    \end{eqnarray}
Note that if we put $f_{u_{i}} =\phi (u_{i})$ then $S=\langle (u_{i},f_{u_{i}})\mid 1\leq i\leq n\rangle$ and one can
easily obtain the associated involutive non-degenerate set theoretical solution $r:X^{2}\rightarrow X^{2}$ defining the
monoid of $I$-type. Indeed, if we set $X=\{ u_{1},\ldots , u_{n}\}$, then  $r(u_{i},u_{j})
=(f_{u_{i}}(u_{j}),f_{f_{u_{i}}(u_{j})}^{-1}(u_{j}))$. Obviously, $\phi (\Fa_{n}) =\langle \phi (a) \mid a \in
\FaM_{n}\rangle =\langle f_{u_{i}} \mid 1\leq i\leq n\rangle$. Note that, because of Proposition~2.2 in \cite{ESS},
$T\inv f_x\inv T=g_x$, where $T:X\rightarrow X$ is the bijective map defined by $T(y)=f_y\inv(y)$. Hence $\GEN{f_x:x\in
X}$ is isomorphic with $\GEN{g_x:x\in X}$.

So, in order to describe all involutive
non-degenerate set theoretical solutions of the
Yang-Baxter equation one needs to characterize the
groups of $I$-type. An important first step in this
direction is to classify the finite groups that are
of the type $\phi (\Fa_{n})$ for some group of
$I$-type ${\cal G}$, as in (\ref{Itypegroup}). A
finite group with this property we will call an {\em
involutive Yang-Baxter} (IYB, for short) group. A
second step is to describe all groups of $I$-type
that have a fixed associated IYB group $G$.

In \cite{ESS}, Etingof, Schedler and Soloviev proved
that any group of $I$-type is solvable. As a
consequence, every IYB group is solvable. In
\cite{ESS} it is also proved that a group ${\cal G}$
is of $I$-type if and only if there is a bijective
1-cocycle ${\cal G}\rightarrow \Fa_n$ with respect
to some action of ${\cal G}$ on $\Fa_n$ which
factors through the natural action of $\Sym_n$ on
$\{u_1,\dots,u_n\}$.

Now, if ${\cal G}=\{ (a,\phi (a)) \mid a\in \Fa_{n}\}$ is a group of $I$-type then the IYB group $G=\phi (\Fa_{n})$
naturally acts   on the quotient group $A=\Fa_{n}/K$, where $K=\{ a\in \Fa_{n} \mid \phi (a)=1\}$ and we obtain a
bijective associated $1$-cocycle $G\rightarrow A$ with respect to this action. By a result of Etingof and Gelaki
\cite{EG}, this bijective 1-cocycle yields a non-degenerate 2-cocycle on the semidirect product $H=A\rtimes G$. This
has been generalized by Ben David and Ginosar \cite{BDG} to more general extensions $H$ of $A$ by $G$
 with a
bijective 1-cocycle from $G$ to $A$. This
construction of Etingof and Gelaki and of Ben David
and Ginosar gives rise to a group of central type in
the sense of \cite{BDG},  i.e. a finite group $H$
with a 2-cocycle $c\in Z^2(H,\C^*)$ such that the
twisted group algebra $\C^cH$ is isomorphic to a
full matrix algebra over the complex numbers, or
equivalently $H=K/Z(K)$ for a finite group $K$ with
an irreducible character of degree
$\sqrt{[K:Z(K)]}$. This provides a nice connection
between IYB groups and groups of central type that
should be investigated. The authors thank Eli
Aljadeff for pointing out this connection.

It is worth  mentioning that the  semigroup algebra
$F S$ of a monoid of $I$-type $S$ over an arbitrary
field $F$ shares many properties with  the
polynomial algebra in finitely many commuting
variables. For example, in \cite{GIVdB}, it is shown
that $FS$ is a domain that satisfies a polynomial
identity and that it is a maximal order in its
classical ring of quotients. In particular, the
group of $I$-type $SS^{-1}$ is finitely generated
abelian-by-finite and torsion free (i.e., it is a
Bieberbach group). The homological properties for
$FS$ were the main reasons for studying monoids of
$I$-type in \cite{GIVdB} and it was inspired by
earlier work of Tate and Van den Bergh on Sklyanin
algebras.

In this paper we investigate group theoretical properties of IYB groups. The  content of the paper is  as follows. In
Section~\ref{SecEquivCond} we obtain several characterizations of IYB groups. These allow us in Section~\ref{ClassIYB}
to prove that the class of IYB groups includes the following: finite abelian-by-cyclic groups, finite nilpotent groups
of class $2$, direct products and wreath products of IYB groups, semidirect products $A\rtimes H$ with $A$ a finite
abelian group $A$ and $H$  an IYB group,  Hall subgroups of IYB groups, Sylow subgroups of symmetric groups $\Sym_n$.
These results imply that any finite nilpotent group is a subgroup of an IYB (nilpotent) group.   It is unclear whether
the class of IYB groups is closed for taking subgroups. As a consequence, we do not know whether the class of IYB
groups contains all finite nilpotent groups. At this point it nevertheless is tempting to conjecture that the class of
IYB groups coincides with that of all solvable finite groups. To prove this, one would like to be able to lift the IYB
structure from subgroups $H$ or quotient groups $\overline{G}$ of a given group $G$ to  $G$. In Section~\ref{ExSect},
we give some examples  of IYB groups of minimal order that are not covered by the general results of
Section~\ref{ClassIYB}.   The last example, a 3-group of class 3, shows that not every IYB homomorphism (see
Section~\ref{SecEquivCond} for the definition) of a quotient of $G$ can be lifted to an IYB homomorphism of $G$. This
indicates that there is no obvious inductive process to prove that nilpotent finite groups are IYB. In
Section~\ref{extsect} we consider the connection between set theoretic solutions of the Yang-Baxter equation and IYB
groups from a different perspective. If $r(x_1,x_2)=(f_{x_1}(x_2),g_{x_2}(x_1))$ is an involutive non-degenerate
solution on a finite set $X$ of the set-theoretical Yang-Baxter equation then it is easy to produce,  in an obvious
manner, infinitely many solutions with the same associated IYB group, namely for every set $Y$ let $r_Y:(X\cup
Y)^2\rightarrow (X\cup Y)^2$ be given by $r_Y((x_1,y_1),(x_2,y_2))=((f_{x_1}(x_2),y_1),(g_{x_2}(x_1),y_2))$. We show an
alternative way of obtaining another IYB map on $X\times X$ with the same associated IYB group. Hence providing, in a
non-obvious fashion, infinitely many set theoretic solutions for the same IYB group.

\section{A characterization of IYB groups}\label{SecEquivCond}

In this section we  obtain several characterizations
of IYB groups. In order to state these we introduce
the following terminology.  For a finite set $X$ we
denote by $\Sym_{X}$ the symmetric group on $X$. An
{\em involutive Yang-Baxter map} (IYB map, for
short) on a finite set $X$ is a map
$\lambda:X\rightarrow \Sym_X$ satisfying
    \begin{equation}\label{InvoYB}
    \lambda(x) \lambda(\lambda(x)\inv(y)) = \lambda(y) \lambda(\lambda(y)\inv(x)) \quad (x,y\in X).
    \end{equation}
The justification for this terminology is based on
the fact that each IYB map yields an involutive non
degenerate set theoretical solution of the
Yang-Baxter equation and conversely. Indeed, let
$r:X^{2}\rightarrow X^{2}$ be a bijective map. As
before, denote $r(x,y)=(f_{x}(y),g_{y}(x))$. From
the proof of \cite[Theorem~4.1]{CJO}, it follows
that $r:X^{2}\rightarrow X^{2}$ is an involutive non
degenerate set theoretical solution of the
Yang-Baxter equation if and only if $f_x\in \Sym_X$
for all $x\in X$ and the map $\lambda\colon
X\rightarrow \Sym_X$ defined by $\lambda (x)=f_x$,
for all $x\in X$, is an IYB map.

\begin{theorem}\label{EquivCond}
The following conditions are equivalent for a finite
group $G$.
\begin{enumerate}
\item
$G$ is an IYB group, that is, there is a map
$\phi:\Fa_n\rightarrow \Sym_n$ such that
$\{(a,\phi(a)):a\in \Fa_n\}$ is a subgroup of
$\Fa_n\rtimes \Sym_n$ and $G$ is isomorphic to
$\phi(\Fa_n)$.
\item
There is an abelian group $A$, an action of $G$ on
$A$ and a group homomorphism $\rho:G\rightarrow
A\rtimes G$ such that $\pi_G\rho=\id_G$ and
$\pi_A\rho:G\rightarrow A$ is bijective, where
$\pi_G$ and $\pi_A$ are the natural projections on
$G$ and $A$ respectively.
\item There is an abelian group $A$, an action of $G$ on $A$ and a bijective 1-cocycle
    $G\rightarrow A$.
\item\label{StrongIYB}
There exists an IYB map $\lambda\colon A\cup
X\rightarrow \Sym_{A\cup X}$ satisfying the
following conditions:
    \begin{enumerate}
    \item $\lambda(A)$ is a subgroup of $\Sym_{A\cup X}$ isomorphic to $G$,
    \item $A\cap X=\emptyset$,
    \item $\lambda(x)=\id_{A\cup X}$ for all $x\in X$,
    \item $\lambda(a)(b)\in A$ for all $a,b\in A$ and
    \item $\lambda|_A$ is injective.
    \end{enumerate}
\item
$G\cong \lambda(X)$ for some IYB map
$\lambda:X\rightarrow \Sym_X$ whose image is a
subgroup of $\Sym_X$.
\item
$G\cong \GEN{\lambda(X)}$ for some IYB map
$\lambda:X\rightarrow \Sym_X$.
\item
There exist a group homomorphism $\mu:G\rightarrow
\Sym_G$ satisfying
    \begin{equation}\label{IYBS}
    x\mu(x)\inv(y)=y\mu(y)\inv(x),
    \end{equation}
for all $x,y\in G$.
\item
There exist a generating subset $Z$ of $G$ and a
group homomorphism $\mu:G\rightarrow \Sym_Z$
satisfying $(\ref{IYBS})$ for all $x,y\in Z$.
\end{enumerate}
\end{theorem}

\begin{proof}
{\it 1 implies 2}.  Suppose that $G$ is an IYB
group. Thus there is a map $\phi\colon
\Fa_n\rightarrow \Sym_n$ such that ${\cal G}=\{ (a,
\phi (a))\mid a\in\Fa_n\}$ is a subgroup of
$\Fa_n\rtimes \Sym_n$, and  $G\cong \phi(\Fa_n)$. We
may assume that $G=\phi(\Fa_n)$. Let $K=\{ (a, \phi
(a))\in {\cal G}\mid\phi(a)=1\}$. Then $G$ is
isomorphic to ${\cal G}/K$ and $B=\{ b\in \Fa_n\mid
(b,1)\in K\}$ is a subgroup of $\Fa_n$. Furthermore,
if $b\in B$ and $a\in \Fa_n$ then
    $$(a,\phi(a))(b,1)(a,\phi(a))\inv =
    (a\phi(a)(b),\phi(a))(\phi(a)\inv(a\inv),\phi(a)\inv) =
    (\phi(a)(b),1)\in K.$$
This shows that $B$ is invariant by the restriction
to $G$ of the action of $\Sym_n$ on $\Fa_n$. Hence
this action induces an action of $G$ on $A=\Fa_n/B$.
Furthermore, $(b,1)(a,\phi(a))=(ba,\phi(a))$ and so
$\phi(ba)=\phi(a)$ for each $b\in B$ and $a\in
\Fa_n$. Let $a,c\in \Fa_n$ such that
$\phi(a)=\phi(c)$. Then
    $$\phi(a)=\phi(cc^{-1}a)=\phi(c)\phi(\phi(c)^{-1}(c^{-1}a)).$$
Hence $\phi(\phi(c)^{-1}(c^{-1}a))=1$, that is
$\phi(c)^{-1}(c^{-1}a)\in B$. Since $B$ is invariant
under the action of the elements of $G$, we have
$c\inv a \in B$, i.e. $\phi(c^{-1}a)=1$. Thus $\phi$
induces a bijective map $\lambda\colon A \rightarrow
G$. Furthermore the map $\rho\colon G\rightarrow
A\rtimes G$ defined by $\rho
(g)=(\lambda^{-1}(g),g)$ satisfies the required
conditions.

{\it 2 implies 3}. Let $G$, $A$ and
$\rho:G\rightarrow A\rtimes G$ as in {\it 2}. Thus
$\rho(G)=\{(\pi(g),g):g\in G\}$ is a subgroup of
$A\rtimes G$, with $\pi=\pi_A\rho$. Then
$(\pi(gh),gh)=(\pi(g),g)(\pi(h),h)=(\pi(g)g(\pi(h)),gh)$
and hence $\pi$ is a 1-cocycle and it is bijective
by assumption.

{\it 3 implies 4}. Let $G$ act on an abelian group
$A$ and let $\pi:G\rightarrow A$ be a bijective
1-cocycle.  If $g,h\in G$ then
    $$g\pi\inv(g\inv \pi(h))=\pi\inv \pi(g\pi\inv(g\inv \pi(h))) =
    \pi\inv(\pi(g) g \pi \pi\inv g\inv \pi(h)) = \pi\inv(\pi(g)\pi(h)).$$
Therefore
    \begin{equation}\label{piinv}
    g\mu(g)\inv(h)=h\mu(h)\inv(g), \quad (g,h\in G),
    \end{equation}
because $A$ is abelian, where $\mu(g)\colon
G\rightarrow G$ is the map given by
$\mu(g)(h)=\pi^{-1}(g\pi(h))$, for $g,h\in G$.

Let $\psi\colon G\rightarrow \Sym_X$ be a
monomorphism, where $X$ is a finite set such that
$A\cap X=\emptyset$. Let $\nu\colon A\cup
X\rightarrow G$ be the map defined by
    $$\nu(t)=\left\{\begin{array}{ll}
    \pi\inv(t), &\mbox{ if }t\in A;\\
    1, &\mbox{ if }t\in X.
     \end{array} \right.$$
Let $\varphi\colon G\rightarrow \Sym_{A\cup X}$ be
the map such that $\varphi(g)$ acts as $g$ on $A$
and as $\psi(g)$ on $X$. It is clear that $\varphi$
is a monomorphism. Clearly $\lambda=\varphi\nu$
satisfies conditions {\it (a)-(e)} and we shall see
that $\lambda$ is an IYB map, that is, it satisfies
condition (\ref{InvoYB}). If $x\in X$ then
$\lambda(x)=\lambda(\lambda(y)^{-1}(x))=1$ and thus
the equality in (\ref{InvoYB}) holds for every $y\in
A\cup X$. Similarly the equality holds if $y\in X$.
Finally, assume that $x,y\in A$. In this case
$\pi\inv(x)\; \pi\inv(\pi\inv(x)^{-1}(y)) =
\pi\inv(y) \; \pi\inv(\pi\inv(y)^{-1} (x))$, by
(\ref{piinv}), and the equality (\ref{InvoYB})
follows by the following calculation
    $$(\lambda(x)\; \lambda(\lambda(x)^{-1}(y)))(v)=
    \left\{\begin{array}{ll}
    \pi\inv(x)\; \pi\inv(\pi\inv(x)^{-1}(y))(v),&\mbox{ if }v\in A;\\
    \psi(\pi\inv(x)\; \pi\inv(\pi\inv(x)^{-1}(y)))(v),&\mbox{ if }v\in X.
    \end{array}\right.$$

{\it 4 implies 5} is obvious.

{\it 5 implies 7}. Suppose that there exists an IYB
map $\lambda\colon X\rightarrow \Sym_X$ such that
$\lambda (X)$ is a subgroup of $\Sym_X$ and $G\cong
\lambda(X)$. We may assume that $G=\lambda (X)$.
Define  $\mu\colon G\rightarrow \Sym_G$ by
$$\mu(g)(\lambda (x))=\lambda (g(x)),$$ for all
$g\in G$ and all $x\in X$. We shall see that $\mu$
is well defined. Let $x,y\in X$ such that $\lambda
(x)=\lambda (y)$ and let $g\in G$. Since $G=\lambda
(X)$, there exists $z\in X$ such that $g=\lambda
(z)^{-1}$. Since $\lambda$ is an IYB map,
$$\lambda(z) \;
\lambda(\lambda(z)^{-1}(x))=\lambda(x)\;
\lambda(\lambda(x)^{-1}(z)) = \lambda(y) \;
\lambda(\lambda(y)^{-1}(z)) = \lambda (z)\;
\lambda(\lambda(z)^{-1}(y)).$$ So
$\lambda(\lambda(z)^{-1}(x))=\lambda(\lambda(z)^{-1}(y))$,
that is $\lambda(g(x))=\lambda(g(y))$. Hence $\mu$
is well defined.

Let $g,h\in G$ and $x\in X$. We have
$$\mu(gh)(\lambda(x))=\lambda
(gh(x))=\mu(g)(\lambda(h(x)))=\mu(g)(\mu(h)(\lambda(x))).$$
Therefore $\mu$ is a group homomorphism.
Furthermore, there exist $y,z\in X$ such that
$g=\lambda(y)$ and $h=\lambda(z)$, and we have
    $$g\; \mu(g)^{-1}(h)=\lambda(y)\; \lambda(\lambda(y)^{-1}(z))=
    \lambda(z)\; \lambda(\lambda(z)^{-1}(y))=h\; \mu(h)^{-1}(g).$$

{\it 7 implies 8} is obvious.

{\it 8 implies 6}. Let $\mu:G=\GEN{Z}\rightarrow
\Sym_Z$ be a map satisfying condition {\it 8}. Let
$\alpha:G\rightarrow \Sym_Y$ be a monomorphism, for
some finite set $Y$ such that $Y\cap Z=\emptyset$.
Let $X=Y\cup Z$. Let $\lambda\colon X\rightarrow
\Sym_X$ be the map defined by
$$\lambda(y)(x)=x,\quad
\lambda(z)(y)=\alpha(z)(y)\quad\mbox{and}\quad\lambda(z)(z')=\mu(z)(z'),$$
for all $x\in X$, $y\in Y$ and $z,z'\in Z$.

Since $\alpha$ and $\mu$ are group homomorphisms,
the restriction of $\lambda$ to $Z$,
$\lambda|_Z\colon Z\rightarrow \Sym_X$, extends to a
group homomorphism $f:G\rightarrow \Sym_X$ which is
injective, because so is $\alpha$. Thus $G\cong
f(G)=\langle\lambda(X)\rangle$, and so it is enough
to show that $\lambda$ is an IYB map, i.e. condition
(\ref{InvoYB}) holds for every $x,y\in X$. If
$x,y\in Y$ then
    $$\lambda(x)=\lambda(y)=\lambda(\lambda(x)^{-1}(y))=\lambda(\lambda(y)^{-1}(x))=\id_X$$
and condition (\ref{InvoYB}) follows easily in this
case. If $x\in Y$ and $y\in Z$ then
    $$\lambda(x)=\lambda(\lambda(y)^{-1}(x))=\id_X.$$
Hence
    $$\lambda(x)\; \lambda(\lambda(x)^{-1}(y))=\lambda(y)=\lambda(y)\; \lambda(\lambda(y)^{-1}(x)).$$
Finally, if $x,y\in Z$ then
$\lambda(x)\inv(y)=\mu(x)\inv(y)$ and
$\lambda(y)\inv(x)=\mu(y)\inv(x)$. Thus $x\;
\lambda(x)\inv(y)=x\; \mu(x)\inv(y)=y\;
\mu(y)\inv(x)=y\; \lambda(y)\inv(x)$. Then
    $$\lambda(x)\; \lambda(\lambda(x)\inv(y)) =
    f(x\; \lambda(x)\inv(y)) =f(y\; \lambda(y)\inv(x))= \lambda(y)\; \lambda(\lambda(y)\inv(x)),$$
as wanted.

{\it 6 implies 1}. Let $\lambda:X\rightarrow \Sym_X$
be an IYB map such that $G\cong \GEN{\lambda(X)}$.
Let $\FaM_X$ be the free abelian monoid with basis
$X$. We extend $\lambda$ to a map
$\lambda:\FaM_X\rightarrow \Sym_X$ by setting
    $$\lambda(1)=1 \quad \text{and} \quad \lambda(xa)=\lambda(x)\;
    \lambda\left(\lambda(x)\inv(a)\right), \text{ if } x\in X
    \text{ and } a\in \FaM_X,$$
where the action of $\Sym_X$ on $\FaM_X$ is the
natural one. We need to show that $\lambda$ is well
defined or equivalently that if $a=xb=yc$ with
$x,y\in X$ and $b,c\in \FaM_X$ then $\lambda(x)\;
\lambda\left(\lambda(x)\inv(b)\right)= \lambda(y)\;
\lambda \left( \lambda(y)\inv(c) \right)$. We may
assume that $x\neq y$. In that case there is $d\in
\FaM_X$ such that $b=yd$ and $c=xd$.  We argue on
$l(d)$, where the map $l:\FaM_X\rightarrow \N$ is
defined by $l(\prod_{x\in X} x^{m_x}) = \sum_{x\in
X} m_x$, ($m_x\ge 0$). Since $\lambda$ is an IYB
map, the claim follows if $l(d)=0$. In the induction
argument we use that $\lambda(a_1)$ is well defined
if $l(a_1)<l(a)$.
    \begin{eqnarray*}
    \lambda(x) \; \lambda\left(\lambda(x)\inv(b)\right)
    &=& \lambda(x) \; \lambda\left(\lambda(x)\inv(y) \; \lambda(x)\inv(d)\right) \hspace{1cm} (\text{because } \lambda(x) \text{ is a
    isomomorphism}) \\
    &=& \lambda(x) \; \lambda(\lambda(x)\inv(y))\;
    \lambda\left( \lambda \left( \lambda(x)\inv(y) \right)\inv \left( \lambda(x)\inv(d) \right) \right)
    \hspace{0.7cm} (\text{by induction}) \\
    &=& \lambda(x) \; \lambda(\lambda(x)\inv(y))\;
    \lambda\left( \left[ \lambda(x) \; \lambda \left( \lambda(x)\inv(y) \right)\right]\inv (d) \right) \\
    &=& \lambda(y) \; \lambda(\lambda(y)\inv(x))\;
    \lambda\left( \left[ \lambda(y) \; \lambda \left( \lambda(y)\inv(x) \right)\right]\inv (d) \right)
    \hspace{1cm} (\text{by hypothesis}) \\
    &=& \lambda(y) \; \lambda(\lambda(y)\inv(x))\;
    \lambda\left( \lambda \left( \lambda(y)\inv(x) \right)\inv \left(\lambda(y)\inv (d)\right) \right) \\
    &=& \lambda(y) \; \lambda\left( \lambda(y)\inv(x) \; \lambda(y)\inv (d)\right)
    \hspace{4.5cm} (\text{by induction}) \\
    &=& \lambda(y) \; \lambda\left( \lambda(y)\inv(c)\right)
    \hspace{2.8cm} (\text{because } \lambda(y) \text{ is a isomomorphism}). \\
    \end{eqnarray*}
So $\lambda$ is well defined.

Now we show that $\lambda(ab) = \lambda(a)
\lambda(\lambda(a)\inv(b))$, for every $a,b\in
\FaM_X$. Again we argue by induction on $l(a)$. The
case $l(a)\le 1$ follows by the definition of
$\lambda$ and the induction step follows by the
following computation,  with $x\in X$:
    \begin{eqnarray*}
    \lambda(xa\cdot b) &=& \lambda(x)\lambda\left(\lambda(x)\inv(ab)\right)
    = \lambda(x)\lambda\left(\lambda(x)\inv(a) \lambda(x)\inv(b)\right) \\
    &=& \lambda(x)\lambda\left(\lambda(x)\inv(a)\right)
    \lambda\left( \lambda\left(\lambda(x)\inv(a)\right)\inv \left(\lambda(x)\inv(b)\right) \right) \\
    &=& \lambda(x)\lambda\left(\lambda(x)\inv(a)\right)
    \lambda\left( \left[\lambda(x) \lambda\left(\lambda(x)\inv(a)\right)\right]\inv (b)\right) \\
    &=& \lambda(xa) \lambda\left( \lambda(xa)\inv (b) \right). \\
    \end{eqnarray*}
By \cite[Lemma~8.2.2]{JObook}, $\lambda$ extends
uniquely to a map $\lambda\colon\Fa_X\rightarrow
\Sym_X$ such that
     $$\lambda(1)=1 \quad \text{and} \quad \lambda(ab) = \lambda(a) \lambda(\lambda(a)\inv(b))$$
for every $a,b\in \Fa_X$,  where $\Fa_X$ is the free
abelian group with basis $X$. Now it is easy to see
that $\{(a,\lambda(a))\mid a\in \Fa_X \}$ is a
subgroup of $\Fa_X\rtimes \Sym_X$. Thus
$\lambda(\Fa_X)$ is an IYB group. By the definition
of $\lambda$ one has
$\GEN{\lambda(X)}=\lambda(\Fa_X)$ and so $G$ is an
IYB group.
\end{proof}

Let $G$ be a finite group. A group homomorphism
$\mu:G\rightarrow \Sym_G$ satisfying condition
$(\ref{IYBS})$ we call an {\em IYB morphism} of $G$.
So, because of Theorem~\ref{EquivCond}, a finite
group is an IYB group if and only if it admits an
IYB morphism, or equivalently, if it admits a
bijective 1-cocycle $G\rightarrow A$ with respect to
some action of $G$ on an abelian group $A$. In the
remainder of the paper these two conditions  turn
out to be very useful.  It therefore is convenient
to know how to pass directly from one to the other.

If $\pi:G\rightarrow A$ is a bijective 1-cocycle
then the equality in (\ref{piinv}) holds for every
$g,h\in G$. Equivalently, for every $g,h\in G$,
 $$g\; (\pi\inv \circ \alpha_g\inv \circ \pi)(h) = h
\; (\pi\inv \circ\alpha_h\inv \circ \pi)(g),$$ where
$\alpha_g\colon A\rightarrow A$ is the map given by
$\alpha_g(a)=g(a)$, for $a\in A$. Thus $g\mapsto
\pi\inv \circ \alpha_g \circ \pi$ is an IYB morphism
of $G$. (Notice that the proof of (\ref{piinv}) does
not use that $\alpha_g$ is a group homomorphism.
Thus even if $\alpha_g$ is not assumed to be a group
homomorphism for all $g$, still this defines an IYB
morphism.)

Conversely, let $\mu$ be an IYB morphism of $G$.  We
define a new product $*$ on $G$. So let $a,b\in G$.
Put
    $$a*b = a\mu(a)\inv(b)=b\mu(b)\inv(a).$$
We claim that $A=(G,*)$ is an abelian group, $\mu$
gives an action of $G$ on $A$ and the identity
$1:G\rightarrow A$ is a bijective 1-cocycle. That
$\mu:G\rightarrow \Sym_A$ is a group homomorphism
and $1(gh)=1(g)\mu(g)(1(h))$ is clear. To prove the
associativity of $*$ we first prove the following
equality:
    \begin{equation}\label{2var}
    \mu(x)(yz)=\mu(x)(y)\cdot \mu(\mu(y)\inv(x\inv))\inv(z) \quad \quad (x,y,z\in G).
    \end{equation}
Indeed
    $$\matriz{{rcl}
    \mu(x)(yz) & = & xx\inv \mu(x\inv)\inv(yz) = x y z \mu(yz)\inv(x\inv) = x y z \mu(z)\inv(\mu(y)\inv(x\inv)) \\
    &=& x y \mu(y)\inv(x\inv) \mu(\mu(y)\inv(x\inv))\inv(z) = x x\inv \mu(x\inv)\inv(y)\mu(\mu(y)\inv(x\inv))\inv(z) \\
    &=& \mu(x)(y)\mu(\mu(y)\inv(x\inv))\inv(z)}.$$
Now, if $x,y,z\in A$ then
    $$\matriz{{rcl}
    x* (y * z) &=& x\mu(x\inv)(y\mu(y\inv)(z)) = x\mu(x\inv)(y)\mu(\mu(y\inv)(x))\inv(\mu(y\inv)(z))\quad\mbox{(by (\ref{2var}))} \\
    &=& x\mu(x\inv)(y)\mu((y\mu(y\inv)(x))\inv(z) = (x* y)\mu(x* y)\inv(z)=(x* y)* z.
    }$$
and the associativity of $*$ follows. The
verification of the remaining axioms of group are
straightforward.

In this paper (as customary), it is assumed that all
actions of a group $G$ on an abelian group $A$ are
via automorphisms on $A$. Hence, we should now
verify that $\mu (g)$ is an automorphism of $A$, for
each $g\in G$. However, because of the presence of
the bijective $1$-cocycle it turns out that this
property automatically is satisfied. More precisely
we prove the following result.

\begin{proposition}\label{Action}
Let $G$ and $A$ be two groups with $A$ abelian. Let
$\alpha:G\rightarrow \Sym_A$  be a group
homomorphism and assume that there is a bijection
$\pi:G\rightarrow A$ satisfying
    \begin{equation}\label{cocycle}
    \pi(gh)=\pi(g)\alpha(g)(\pi(h)) \quad \text{for every } g,h\in G.
    \end{equation}
Then $\alpha$ defines an action (by automorphisms)
of $G$ on $A$.
\end{proposition}

\begin{proof}
Let $\alpha_g$ denote the image of $g\in G$ by
$\alpha$. We have to show that $\alpha_g$ is a group
homomorphism for every $g\in G$. We define
$\mu:G\rightarrow \Sym_G$ by $\mu(g)=\pi\inv
\alpha_g \pi$. As it was pointed above, even if
$\alpha_g$ is not assumed to be a group
homomorphism, $\mu$ is an IYB morphism. In
particular $\mu$ satisfies (\ref{2var}). Then
    $$\matriz{{rcll}
    \alpha_g(ab) &=&
    \alpha_g(\pi(\pi\inv(a)) \alpha_{\pi\inv(a)}\pi\pi\inv\alpha_{\pi\inv(a)}\inv(b)) \\
    &=&\alpha_g(\pi(\pi\inv(a) \pi\inv\alpha_{\pi\inv(a)}\inv(b))) & (\ref{cocycle})\\
    &=& \pi \mu(g)(\pi\inv(a) \pi\inv\alpha_{\pi\inv(a)}\inv(b)) & (\text{definition of } \mu)\\
    &=& \pi (\mu(g)(\pi\inv(a)) \mu(\mu(\pi\inv(a))\inv(g\inv))\inv \pi\inv\alpha_{\pi\inv(a)}\inv(b))
            & (\ref{2var})\\
    &=& \pi (\mu(g)(\pi\inv(a)) \pi\inv \alpha_{\mu(\pi\inv(a))\inv(g\inv)}\inv \alpha_{\pi\inv(a)}\inv(b))
            & (\text{definition of } \mu)\\
    &=& \pi (\mu(g)(\pi\inv(a)) \pi\inv \alpha_{\pi\inv(a)\mu(\pi\inv(a))\inv(g\inv)}\inv (b))
            & (\alpha \text{ is homomorphism}) \\
    &=& \pi (\mu(g)(\pi\inv(a)) \pi\inv \alpha_{g\inv \mu(g)(\pi\inv(a))}\inv (b))
            & (\mu \text{ is IYB morphism})\\
    &=& \pi (\mu(g)(\pi\inv(a))) \alpha_{\mu(g)(\pi\inv(a))} \pi \pi\inv \alpha_{g\inv \mu(g)(\pi\inv(a))}\inv (b)
            & (\ref{cocycle}) \\
    &=& \alpha_g(a) \alpha_g (b) & (\alpha \text{ is homomorphism})
    }$$
\end{proof}

\section{The class of IYB groups}
\label{ClassIYB}

In this section we prove that some classes of groups
 consist of  IYB groups and that the class of
IYB groups is closed under some constructions. We
start with some easy consequences of the
characterization of IYB groups in terms of
1-cocycles.

 Recall from \cite[Theorem 2.15]{ESS} that
every group of $I$-type is solvable, and henceforth
so is every IYB group. For the sake of completeness,
we include a proof, which is essentially  the proof
of \cite[Theorem 2.15]{ESS}. Indeed, if
$\pi:G\rightarrow A$ is a bijective 1-cocycle then
the 1-cocycle condition implies that if $B$ is a
characteristic subgroup of $A$ then $\pi\inv(B)$ is
a subgroup of $G$.  In particular, if $p$ is a prime
and $P$ is the Hall $p'$-subgroup of $A$, then
$\pi\inv(P)$ is a Hall $p'$-subgroup of $G$. By a
theorem of P. Hall \cite[9.1.8]{Robinson}, $G$ is
solvable.

\begin{corollary}\label{Hall}
If $G$ is an IYB group then its Hall subgroups are
also IYB.
\end{corollary}

\begin{proof}
Assume that $G$ is an IYB group. By
Theorem~\ref{EquivCond} there is a bijective
1-cocycle $\pi:G\rightarrow A$, for some abelian
group $A$. If $B$ is a Hall subgroup of $A$ then $B$
is invariant under the action of $G$ and hence
$H=\pi\inv(B)$ is a subgroup of $G$ of the same
order than $B$.  Then the restriction of $\pi$ to
$H$, $\pi|_H:H\rightarrow B$, is a bijective
1-cocycle, with respect to the action of $H$
restricted to $B$.
\end{proof}

\begin{corollary}\label{directproduct}
The class of IYB groups is closed under direct
products.
\end{corollary}

\begin{proof}
Let $G_1$ and $G_2$ be IYB groups. By
Theorem~\ref{EquivCond} there are bijective
1-cocycles $\pi_i:G_i\rightarrow A_i$ with respect
to some action of $G_i$ on an abelian group $A_i$
($i=1,2$). Then $\pi_1\times \pi_2:G_1\times
G_2\rightarrow A_1\times A_2$ is a bijective
1-cocycle with respect to the obvious action of
$G_1\times G_2$ on $A_1\times A_2$, namely
$(g_1,g_2)(a_1,a_2)=(g_1(a_1),g_2(a_2))$.
\end{proof}

The next two results provide more closure properties
of the class of IYB groups.

\begin{theorem}\label{semidirect}
Let $G$ be a finite group such that $G=AH$, where
$A$ is an abelian normal subgroup of $G$ and $H$ is
an IYB subgroup of $G$. Suppose that there is a
bijective 1-cocycle $\pi:H\rightarrow B$, with
respect to an action of $H$ on the abelian group $B$
such that $H\cap A$ acts trivially on $B$. Then $G$
is an IYB group.

In particular, every semidirect product $A\rtimes H$
of a finite abelian group $A$ by an IYB group $H$ is
IYB.
\end{theorem}

\begin{proof}
Let $N=\{ (h^{-1},\pi(h))\in H\times B \mid h\in
H\cap A\}$. Since $H\cap A$ acts trivially on $B$,
$\pi(ah)=\pi(a)\pi(h)$ for every $a\in A\cap H$ and
$h\in H$. It follows that $N$ is a subgroup of
$A\times B$.

Let $C=(A\times B)/N$ and let $\overline{(a,b)}$
denote the class of $(a,b)\in A\times B$ modulo $N$.
Note that $|G|=|C|$. Define the following action of
$G$ on $C$:
    $$g\overline{(a,b)}=\overline{(gag^{-1},h(b))},$$
for all $a\in A$, $b\in B$ and $g=a'h\in G$, with
$a'\in A$ and $h\in H$. We shall see that this is
well defined and it  indeed is an action. Let
$a,a'\in A$ and $b,b'\in B$ such that
$\overline{(a,b)}=\overline{(a',b')}$.  We have
$h^{-1}=a^{-1}a'\in H\cap A$ and $\pi(h)=b^{-1}b'$.
Let $g=a_1h_1=a_2h_2\in G$, with $a_1,a_2\in A$ and
$h_1,h_2\in H$. Since $A$ is an abelian normal
subgroup of $G$, we have
    $$(gag^{-1})^{-1}ga'g^{-1}=h_1a^{-1}a'h_1^{-1}=h_1h^{-1}h_1^{-1}\in
    H\cap A$$
and
    $$h_2^{-1}h_1=h_2^{-1}(h_1h_2^{-1})h_2=h_2^{-1}(a_1^{-1}a_2)h_2\in H\cap A.$$
Since $h, h_2^{-1}h_1\in H\cap A$ and $H\cap A$ acts
trivially on $B$, we have
\begin{eqnarray*}
\pi(((gag^{-1})^{-1}ga'g^{-1})^{-1})&=&
 \pi(h_1hh_1^{-1})=\pi(h_1)h_1(\pi(hh_1^{-1}))\\
&=&\pi(h_1)h_1(\pi(h)\pi(h_1^{-1})) =
\pi(h_1)h_1(\pi(h))h_1(\pi(h_1^{-1}))\\ &=&
\pi(h_1)h_1(\pi(h_1^{-1}))h_1(\pi(h)) =
\pi(h_1h_1^{-1})h_1(\pi(h))\\ &=&
h_1(\pi(h))=h_1(b^{-1}b')\\
&=&h_1(b)^{-1}h_1(b')=h_2(b)^{-1}h_1(b').
\end{eqnarray*}
Hence
$\overline{(gag^{-1},h_2(b))}=\overline{(ga'g^{-1},h_1(b'))}$.
Furthermore, if $g,g'\in G$, $g=bh$ and $g'=b'h'$,
with $b,b'\in A$ and $h,h'\in H$, then we have
$gg'=(bhb'h^{-1})hh'$ and
$$(gg')\overline{(a,b)}=\overline{((gg')a(gg')^{-1},(hh')(b))}
=\overline{(g(g'ag'^{-1})g^{-1},h(h'(b)))}=
g(g'\overline{(a,b)}).$$ Therefore we have a well
defined action of $G$ on $C$.

We define $\bar \pi\colon G\rightarrow C$ by
    $$\bar \pi(g)=\overline{(a,\pi(h))}$$
for all $g\in G$ with $g=ah$, where $a\in A$ and
$h\in H$. Note that if $g=ah=a'h'$ with $a,a'\in A$
and $h,h'\in H$, then $a^{-1}a'=hh'^{-1}\in H\cap
A$, $h^{-1}h'=h'^{-1}(hh'^{-1})^{-1}h'\in H\cap A$
and hence
    \begin{eqnarray*}\pi((a^{-1}a')^{-1})&=&\pi(h'h^{-1})
    =\pi(h')h'(\pi(h^{-1}))\\
    &=& \pi(h')hh^{-1}h'(\pi(h^{-1}))=
    \pi(h')h(\pi(h^{-1}))\\
    &=& \pi(h')\pi(h)^{-1},
    \end{eqnarray*}
because $H\cap A$ acts trivially on $B$. Hence
$\bar\pi$ is well defined. Let $g_1=a_1h_1$ and
$g_2=a_2h_2$, with $a_1,a_2\in A$ and $h_1,h_2\in
H$. We have
    \begin{eqnarray*}
    \bar{\pi}(g_1g_2) &=& \bar{\pi}((a_1h_1a_2h_1^{-1})h_1h_2) =
    \overline{(a_1h_1a_2h_1^{-1},\pi(h_1h_2))} \\
    &=& \overline{(a_1h_1a_2h_1^{-1},\pi(h_1)h_1(\pi(h_2)))} =
    \overline{(a_1,\pi(h_1))} \; \overline{(h_1a_2h_1^{-1},h_1(\pi(h_2)))} \\
    &=& \overline{(a_1,\pi(h_1))} g_1 \overline{(a_2,\pi(h_2))} = \bar{\pi}(g_1) g_1(\bar{\pi}(g_2)).
    \end{eqnarray*}
Hence $\bar{\pi}$ is a bijective 1-cocycle and we
conclude that $G$ is an IYB group, by
Theorem~\ref{EquivCond}.
\end{proof}

\begin{theorem}\label{CompSemi}
Let $N$ and $H$ be IYB groups and let
$\pi_N:N\rightarrow A$ be a bijective 1-cocycle with
respect to an action of $N$ on an abelian group $A$.
  If $\gamma:H\rightarrow \Aut(N)$ and
$\delta:H\rightarrow \Aut(A)$ are actions of $H$ on
$N$ and $A$ respectively such that $\delta(h)\pi_N =
\pi_N \gamma(h)$ for every $h\in H$, then the
semidirect product $N\rtimes H$, with respect to the
action $\gamma$, is an IYB group.
\end{theorem}

\begin{proof}
Let $\alpha:N\rightarrow \Aut(A)$   be an action of
$N$ on the abelian group $A$ such that
$\pi_N:N\rightarrow A$ is a bijective 1-cocycle. Let
$\pi_H:H\rightarrow B$ be a bijective 1-cocycle with
respect to the action $\beta:H\rightarrow \Aut(B)$.
We are going to denote by $\alpha_n$, $\beta_h$,
$\gamma_h$ and $\delta_h$ to the images of $n\in N$
or $h\in H$  under the respective maps $\alpha$,
$\beta$, $\gamma$ and $\delta$.

If $h\in H$ and $n_1,n_2\in N$ then
    \begin{eqnarray*}
    \delta_h \pi_N(n_1) \; \delta_h \alpha_{n_1}(\pi_N(n_2))
    &=& \delta_h (\pi_N(n_1) \; \alpha_{n_1}(\pi_N(n_2))) = \delta_h \pi_N(n_1 n_2) = \pi_N \gamma_h(n_1 n_2) \\
    &=& \pi_N(\gamma_h(n_1) \gamma_h(n_2) ) = \pi_N \gamma_h(n_1) \alpha_{\gamma_h(n_1)}(\pi_N \gamma_h(n_2)) \\
    &=& \delta_h \pi_N(n_1) \alpha_{\gamma_h(n_1)} \delta_h(\pi_N(n_2)).
    \end{eqnarray*}
This shows that
   \begin{eqnarray}\label{deltaalpha}
   \delta_h \alpha_n &=& \alpha_{\gamma_h(n)} \delta_h
   \end{eqnarray}
for every $h\in H$ and $n\in N$.

Now we define the following map $\sigma:G=N\rtimes H
\rightarrow \Aut(A\times B)$:
    $$\sigma_{nh}(a,b)=(\alpha_n \delta_h(a),\beta_h(b)),$$
for all $n\in N$ and $h\in H$. Since both $\alpha_n$
and $\delta_h$ are automorphisms of $A$ and
$\beta_h$ is an automorphism of $B$, $\sigma_{nh}$
is an automorphism of $A\times B$. Now we check that
$\sigma$ is a group homomorphism. Let $n_1,n_2\in N$
and $h_1,h_2\in H$.  Then
    \begin{eqnarray*}
    \sigma_{n_1 h_1 n_2 h_2}(a,b) & = & \sigma_{n_1 \gamma_{h_1}(n_2) h_1 h_2}(a,b) =
    (\alpha_{n_1\gamma_{h_1}(n_2)} \delta_{h_1h_2}(a),\beta_{h_1h_2}(b)) \\
    & = & (\alpha_{n_1} \alpha_{\gamma_{h_1}(n_2)} \delta_{h_1} \delta_{h_2}(a),\beta_{h_1} \beta_{h_2}(b))
    = (\alpha_{n_1} \delta_{h_1} \alpha_{n_2} \delta_{h_2}(a),\beta_{h_1} \beta_{h_2}(b))\quad
    \mbox{(by (\ref{deltaalpha}))} \\
    & = & \sigma_{n_1h_1} \sigma_{n_2h_2}(a,b).
    \end{eqnarray*}
Thus $\sigma$ is an action of $G$ on $A\times B$.

Let $\pi:G\rightarrow A\times B$ be given by
$\pi(nh)=(\pi_N(n),\pi_H(h))$. Since both
$\pi_N:N\rightarrow A$ and $\pi_H:H\rightarrow B$
are bijective, so is $\pi$. Moreover
    \begin{eqnarray*}
    \pi(n_1h_1 n_2h_2) &=& \pi(n_1 \gamma_{h_1}(n_2) h_1 h_2) = (\pi_N(n_1 \gamma_{h_1}(n_2)), \pi_H(h_1 h_2)) \\
    & = & (\pi_N(n_1) \alpha_{n_1}\pi_N \gamma_{h_1}(n_2), \pi_H(h_1) \beta_{h_1}(\pi_H(h_2)) ) \\
    & = & (\pi_N(n_1) \alpha_{n_1} \delta_{h_1} \pi_N (n_2), \pi_H(h_1) \beta_{h_1}(\pi_H(h_2)) ) \\
    & = & (\pi_N(n_1), \pi_H(h_1)) (\alpha_{n_1} \delta_{h_1} \pi_N (n_2), \beta_{h_1}(\pi_H(h_2)) ) \\
    & = & \pi(n_1h_1) \sigma_{n_1h_1} (\pi(n_2h_2)), \\
    \end{eqnarray*}
for all $n_1,n_2\in N$ and $h_1,h_2\in H$. Thus
$\pi$ is a bijective 1-cocycle and we conclude that
$G$ is an IYB group. \end{proof}

\begin{corollary}\label{wreath}
Let $G$ be an IYB group and $H$ an IYB subgroup of
$\Sym_n$. Then the wreath  product $G\wr H$ of $G$
and $H$ is an IYB group.
\end{corollary}

\begin{proof}
Recall that $W=G\wr H$ is the semidirect product
$G^n \rtimes H$, where the action of $H$ on $G^n$ is
given by
$\gamma_h(g_1,\dots,g_n)=(g_{h(1)},\dots,g_{h(n)})$.
Since $G$ is an IYB group, there is an action on an
abelian group $A$ which admits a bijective 1-cocycle
$\pi:G\rightarrow A$. This action, applied
componentwise, induces an action $\alpha$ of $G^n$
on $A^n$ and the map $\bar{\pi}:G^n \rightarrow
A^n$, which acts as $\pi$ component wise,  is a
bijective 1-cocycle (see the proof of
Corollary~\ref{directproduct}). Furthermore the map
$\delta:H\rightarrow \Aut(A^n)$, given by
$\delta_h(a_1,\dots,a_n)=(a_{h(1)},\dots,a_{h(n)})$,
for all $h\in H$ and $(a_1,\dots,a_n)\in A^n$, is an
action of $H$ on $A^n$ such that $\delta_h\bar\pi
=\bar\pi \gamma_h$ for every $h\in H$. Thus $W$ is
an IYB group, by Theorem~\ref{CompSemi}.
\end{proof}

\begin{corollary}
Let $n$ be a positive integer. Then the Sylow
subgroups of $\Sym_n$ are IYB groups.
\end{corollary}

\begin{proof} It is known that the Sylow $p$-subgroups
of $\Sym_n$ are isomorphic to a group of the form
$G_1\times G_2$ or $G_3\wr C_p$, where $G_1,G_2,G_3$
are Sylow $p$-subgroup of
$\Sym_{m_1},\Sym_{m_2},\Sym_{m_3}$ respectively, for
some $m_1,m_2,m_3<n$, \cite[pages 10,11]{passman}.
Since $C_p$ is an IYB group, the result follows from
Corollaries~\ref{directproduct} and \ref{wreath} by
induction on $n$.
\end{proof}

As a consequence of this result and
Corollary~\ref{directproduct}, we have the following
result.

\begin{corollary} Any finite nilpotent group is isomorphic to a
subgroup of an IYB  (nilpotent) group. $\qed$
\end{corollary}

 The next result yields many examples of IYB
groups.

\begin{theorem}\label{CarIYB}
Let $G$ be a finite group having a normal sequence
    $$1 = G_0 \lhd G_1 \lhd G_2 \lhd \dots \lhd G_{n-1} \lhd G_n = G$$
satisfying the following conditions:
\begin{itemize}
\item[(i)] For every $1\le i\le n$, $G_i = G_{i-1} A_i$ for some abelian
subgroup $A_i$.
\item[(ii)] $(G_{i-1}\cap (A_i \cdots A_n),G_{i-1})=1$.
\item[(iii)] $A_i$ is normalized by $A_j$ for every $i\le j$.
\end{itemize}

Then $G$ is an IYB group.
\end{theorem}

\begin{proof}
Assumption $(iii)$ implies that for every $1\le i,j
\le n$, either $A_i$  normalizes $A_j$ or $A_j$
normalizes $A_i$. Thus $A_iA_j$ is a subgroup of $G$
and $A_iA_j=A_jA_i$. Furthermore $G_i=A_1\dots A_i$
and, in particular, $Z=A_1\cup \dots \cup A_n$ is a
generating subset of $G$. We set $H_i=A_{i+1}\dots
A_n$ for every $1\le i \le n$. Then $A_i$ normalizes
$G_j$ and $H_i$ normalizes $A_j$ for every $j\le i$.
In particular, $G_i\unlhd G$ for every $i$. Every
$g\in G$ can be written, in a non necessarily unique
way, as $g=g_{i1}g_{i2}$ with $g_{i1}\in G_i$ and
$g_{i2}\in H_i$.

Define $\mu\colon G\rightarrow \Sym_Z$ by
    $$\mu(g)(a)=g_{(i-1)2} a g_{(i-1)2}\inv,$$
for every $g\in G$ and $a\in A_i$  , with $i\geq 1$.
We need to show that the map $\mu$ is well defined.
So let $g=a_1\dots a_n=b_1\dots b_n$, with
$a_j,b_j\in A_j$, for every $j=1,\dots,n$ and let
$a\in A_i$. We can take $g_{i1}=a_1\dots a_i \in
G_i$ and $g_{i2}=a_{i+1}\dots a_n\in H_i$, or we can
take $\bar{g}_{i1}=b_1\dots b_i\in G_i$ and
$\bar{g}_{i2}=b_{i+1}\dots b_n\in H_i$. Since
$G_i\unlhd G$ and $g_{i1}g_{i2}=\bar{g}_{i1}
\bar{g}_{i2}$,
    $$\bar{g}_{i2}\inv g_{i2} = \bar{g}_{i2}\inv g_{i1}\inv \bar{g}_{i1} \bar{g}_{i2} \in G_i\cap H_i.$$
Hence, by condition $(ii)$, $\bar{g}_{i2}^{-1}
g_{i2} a g_{i2}^{-1} \bar{g}_{i2} = a$. Thus
    $$g_{i2} a g_{i2}\inv = \bar{g}_{i2} a \bar{g}_{i2}\inv \in A_i.$$
Since $A_i$ is abelian, conjugating by $a_i$ on the
left side of the previous equality and by $b_i$ on
the right side we have $$g_{(i-1)2} a g_{(i-1)2}\inv
= \bar{g}_{(i-1)2} a \bar{g}_{(i-1)2}\inv .$$
Suppose that $a\in A_i\cap A_{i'}$ with $i'< i$.
Since $H_i$ normalizes both $A_i$ and $A_{i'}$ and
$A_i$ is abelian we have
    $$g_{(i-1)2} a {g_{(i-1)2}}^{-1} = a_i g_{i2} a g_{i2}\inv a_i\inv \in A_{i'}\cap A_i\subseteq G_{i-1}\cap A_i.$$
Since $g_{(i'-1)2}
g_{(i-1)2}\inv=a_{i'}a_{i'+1}\cdots a_{i-1} \in
G_{i-1}$, by $(ii)$,
    $$g_{(i-1)2}ag_{(i-1)2}^{-1}=g_{(i'-1)2} a g_{(i'-1)2}^{-1}.$$
Therefore $\mu$ is well defined.

Let $g,h\in G$ and $a\in A_i$. Then, with the above
notation,
    \begin{eqnarray*}
    \mu(gh)(a) &=& \mu(g_{i1}g_{i2}h_{i1}h_{i2})(a) = \mu((g_{i1}g_{i2}h_{i1}g_{i2}\inv) (g_{i2} h_{i2}))(a) \\
    & = & \mu(g_{i2} h_{i2})(a) = g_{i2} h_{i2} a h_{i2}\inv g_{i2}\inv = \mu(g) \mu(h)(a).
    \end{eqnarray*}
Hence $\mu$ is a group homomorphism.

Let $a\in A_i$ and $b\in A_j$ with $i\leq j$. Then
$$a\mu(a)^{-1}(b)=ab=bb^{-1}ab=b\mu(b)^{-1}(a).$$ By
Theorem~\ref{EquivCond}, $G$ is an IYB group.
 \end{proof}

\begin{corollary}\label{Meta}
Let $G$ be a finite group. If $G=NA$, where $N$ and
$A$ are two abelian subgroups of $G$ and $N$ is
normal in $G$, then $G$ is an IYB group.

In particular, every abelian-by-cyclic finite group
is IYB.
\end{corollary}

\begin{corollary}\label{class2}
Every finite nilpotent group of class $2$ is IYB.
\end{corollary}

\begin{proof}
Let $G$ be a finite nilpotent group of class $2$.
Then there exist $x_1,\dots ,x_n\in G$ such that
$G/Z(G)$ is the inner direct product of $\langle
\bar x_1\rangle,\dots, \langle\bar x_n \rangle$,
where $\bar x$ denotes the class of $x\in G$ modulo
its center $Z(G)$. Let $A_i=\langle \{ x_i\}\cup
Z(G)\rangle$, for all $i=1,\dots ,n$. Let
$G_i=A_1\cdots A_i$, for all $i=1,\dots , n$. It is
easy to check that the group $G$ and the subgroups
$A_i$ and $G_i$ satisfy the hypothesis of
Theorem~\ref{CarIYB}. Hence $G$ is an IYB group.
\end{proof}

\section{Examples} \label{ExSect}

In Section~\ref{ClassIYB} we  have given some
sufficient conditions for a finite group to be IYB.
In this section we present some examples of IYB that
are not covered by these results.

\begin{example}\label{Q8C3}
{\rm  A group of smallest order not satisfying the
conditions of Theorem~\ref{CarIYB} is
     $$G=Q_8\rtimes C_3 = \GEN{x,y,a\mid x^4=x^2y^2=a^3=1, x^y=x\inv, axa\inv=y,aya\inv=xy}.$$
One may try to show that $G$ is IYB by using
Theorem~\ref{semidirect} and Corollary~\ref{Meta}.
 However, the straightforward approach does not
work. Nevertheless, we can still show that $G$ is
IYB as follows.

Let $Z=\{x,x\inv,y,y\inv,xy,yx,a\}$ and define
$\mu\colon Z\rightarrow
 \Sym_Z$ by
 \begin{eqnarray*}    &&\mu(x)=\mu(x^{-1})=(x,x^{-1})(y,y\inv),\quad
 \mu(y)=\mu(y\inv)=(y,y^{-1})(xy,yx),\\
&&\mu(xy)=\mu(yx) =(x,x\inv)(xy,yx),\quad
     \mu(a)=(x,y,xy)(x\inv,y\inv,yx).\end{eqnarray*}
Note that $\mu(x)\mu(y)=\mu(xy)$,
$\mu(x)^4=\mu(x)^2\mu(y)^2=\mu(a)^3=\id$,
$\mu(x)^{\mu(y)}=\mu(x\inv),
\mu(a)\mu(x)\mu(a)\inv=\mu(y),
\mu(a)\mu(y)\mu(a)\inv=\mu(x)\mu(y)$. Hence $\mu$
extends to a homomorphism $\mu\colon Q_8\rtimes
C_3\rightarrow \Sym_Z$. It is easy to check that
$$u\mu(u^{-1})(v)=v\mu(v)^{-1}(u)$$ for all $u,v\in
Z$. Therefore, $Q_8\rtimes C_3$ is an IYB group, by
Theorem~\ref{EquivCond}. $\qed$}\end{example}

\begin{example}\label{CbyAb2}{\rm
By Corollary~\ref{Meta}, every abelian-by-cyclic
group is IYB. However it is not clear whether every
cyclic-by-abelian group is IYB. In this example we
show  that some class of cyclic-by-abelian groups
consists of  IYB groups. It includes all
cyclic-by-two generated abelian $p$-groups.

Consider the following cyclic-by-abelian group of
order $nq_1q_2$:
     $$G=\GEN{a,b,c\mid a^n=1,bab\inv = a^{r_1}, cac\inv = a^{r_2}, b^{q_1}=a^{s_1}, c^{q_2}=a^{s_2}, cbc\inv = a^t b}.$$
Then the following conditions  hold, where $o_n(r)$
denotes the multiplicative order of $r$ modulo $n$
(for $(r,n)=1$):
     $$\matriz{{ll}o_n(r_i)|q_i, & (\mbox{because } (a,b_i^{q_i})=1)
     \\ s_i(r_i-1) \equiv 0 \mod n, & (\mbox{because } (a^{s_1},b)=1=(a^{s_2},c)) \\
     t(1+r_1+r_1^2+\dots+r_1^{q_1-1})\equiv s_1(r_2-1) \mod n, &
(\mbox{because } (a^tb)^{q_1} = ca^{s_1}c\inv)\\
     -t(1+r_2+r_2^2+\dots+r_2^{q_2-1})\equiv s_2(r_1-1) \mod n
    &(\mbox{because } (a^{-t}c)^{q_2} = ba^{s_2}b\inv).}$$
We also assume that there is an integer $u$ such
that
     \begin{equation}\label{rr}
     r_1-1\equiv u(r_2-1) \mod n.
     \end{equation}

Let $A=\GEN{a}$. Set
$Z=\{a,a^2,\dots,a^{n-1},b,ab,\dots,a^{n-1}b,c\}$, a
generating subset of $G$, and let $f,g\in \Sym_Z$ be
given by
     $$f:\left\{\matriz{{ll} a^i\mapsto a^{ir_1}, & (1\le i<n) \\ a^jb
\mapsto a^{tu+r_1j}b, & (0\le j<n) \\
     c\mapsto c} \right.$$
     $$g(x)=cxc\inv, \mbox{ for each } x\in Z.$$
We claim that
    \begin{equation}\label{fg}
    f^{q_1}=g^{q_2}=(f,g)=1.
    \end{equation}
Since $c^{q_2}=a^{s_2}$ and $s_2(r_2-1)\equiv 0 \mod
n$ and $s_2(r_1-1)\equiv us_2(r_2-1)\equiv 0 \mod
n$, one has that $c^{q_2}\in Z(G)$ and so
$g^{q_2}=1$. Notice that $f(a^i)=ba^ib\inv$ for each
$i$. Since $b^{q_1}\in \GEN{a}$, one has
$f^{q_1}(a^i)=a^i$. On the other hand
     $$tu(1+r_1+r_1^2+\cdots+r_1^{q_1-1}) \equiv us_1(r_2-1) \equiv s_1(r_1-1) \equiv 0 \mod n.$$
Using this and that $o_n(r_1)|q_1$, one has
     $$f^{q_1}(a^j b) = a^{tu(1+r_1+r_1^2+\cdots+r_1^{q_1-1})+r_1^{q_1}j}b =a^jb.$$
This shows that $f^{q_1}=1$. Now we check that
$gf=fg$. Since the actions of $f$ and $g$ on the
powers of $a$ are by conjugation by $b$ and $c$
respectively, the action of $(f,g)$ on  these powers
is by conjugation by $(b,c)=a^t$ and so
$(f,g)(a^i)=a^i$. Finally
\begin{eqnarray*}
gf(a^j b) &=&
g(a^{tu+r_1j}b)=a^{(tu+r_1j)r_2+t}b=a^{tu(r_2-1)+tu+r_1r_2j+t}b\\
&=&
a^{t(r_1-1)+tu+r_1r_2j+t}b=a^{tu+r_1(r_2j+t)}b=f(a^{r_2j+t}b)=fg(a^jb).
\end{eqnarray*}
This proves the claim.

Since $G/A\cong C_{q_1}\times C_{q_2}$, there is a
group homomorphism $\mu:G\rightarrow \Sym_Z$ such
that $\mu(a)=\id$, $\mu(b)=f$ and $\mu(c)=g$. Now we
check (\ref{IYBS}) for all $x,y\in Z$. If
 $x\in \GEN{a}$, then $\mu(x)=\id$ and $\mu(y)(x)=yxy\inv$. Thus $x
\mu(x)\inv(y)=xy=y(y\inv x y) = y\mu(y)\inv(x)$ as
wanted.
 If $x=c$ then $\mu(y)(x)=x$ and $\mu(x)(y)=xyx\inv$. Again
$x\mu(x)\inv(y)=xx\inv y x = yx=y\mu(y)\inv(x)$. By
symmetry, (\ref{IYBS}) also holds if $y\in
\GEN{a}\cup \{c\}$. If $v$ is an inverse of $r_1$
modulo $n$ then
 $f\inv(a^jb)=a^{v(j-tu)}b$. Thus, for $x=a^ib,y=a^jb$, one has
     $$\matriz{{cc}
     x\mu(x)\inv(y)=xf\inv(y)=a^ib a^{v(j-tu)} b = a^{i+r_1 v(j-tu)} b^2 =
a^{i+j-tu} b^2 \\
     = a^{j+r_1 v(i-tu)} b^2 = a^jb a^{v(i-tu)} b = y f\inv(x) = y
\mu(y)\inv(x).}$$

We conclude that if condition (\ref{rr}) holds then,
by Theorem~\ref{EquivCond}, $G$ is an IYB group.
Notice that if $n$ is a prime power then, by
interchanging the roles of $b$ and $c$ if needed,
one may assume that condition (\ref{rr}) holds
because the lattice of additive subgroups of $\Z_n$
is linearly ordered. So if $G$ has a normal cyclic
$p$-subgroup $A$ such that $G/A$ is $2$-generated
and abelian then $G$ is IYB. $\qed$}\end{example}

\begin{example}\label{C35}{\rm
By Corollaries~\ref{Meta} and \ref{class2} every
nilpotent group of class at most $2$ is IYB and
every group which is a product of two abelian
subgroups, one of them normal, is IYB. In this
example we consider  3-group of minimal order not
satisfying any of these properties.  Namely let $G$
be the group given by the following presentation
    $$G=\GEN{a,b,c,d,e\mid a,b\in Z(G), dc=acd, ec=bce, ed = cde, a^3=b^3=c^3=d^3=e^3=1}.$$
If $i_1,i_2\in \{ 0,1,2\}, j_1,j_2\in \{ 0,1\}$,
then
\begin{eqnarray}\label{diej}
d^{i_1}e^{j_1}d^{i_2}e^{j_2} =
a^{2i_2(i_2-1)j_1+i_1i_2j_1} c^{i_2 j_1} d^{i_1+i_2}
e^{j_1+j_2}.
\end{eqnarray}

Let $H=\GEN{G',d}$ and $Z=H\cup He$. Let $D,E\in
S_Z$ be given by
    $$D(nd^ie^j)=dnd\inv a^{i+2ij} b^{2i} c^{2j} d^i e^j \quad \text{and} \quad
    E(nd^ie^j)=ene\inv d^{i+j} e^j, \quad (n\in G',i=0,1,2; j=0,1).$$
If $n\in G'$, $i=0,1,2$ and $j=0,1$ then
\begin{eqnarray}\label{Dk}
D^k(nd^ie^j)=D^k(n)D^k(d^i)D^k(e^j)=d^knd^{-k}
a^{k(i+2ij)+k(k-1)j} b^{2ik} c^{2kj} d^i e^j
\end{eqnarray}
and
\begin{eqnarray}\label{Ek}
E^k(nd^ie^j) = E^k(n) E^k(d^i) E^k(e^j) = e^k n
e^{-k} d^{i+kj} e^j,
\end{eqnarray}
for all positive integer $k$. Moreover,
    $$\matriz{{rcl}
    ED(nd^ie^j) & = & E(dnd\inv a^{i+2ij} b^{2i} c^{2j} d^i e^j) =
    ednd\inv e\inv a^{i+2ij} b^{2(i+j)} c^{2j} d^{i+j} e^j \\
    &= & c(dene\inv d\inv)c\inv a^{i+2ij} b^{2(i+j)} c^{2j} d^{i+j} e^j =
    dene\inv d\inv a^{i+2ij} b^{2(i+j)} c^{2j}d^{i+j} e^j  \\
    &= & D(ene\inv d^{i+j} e^j) = DE(nd^ie^j)}.
    $$
Thus $ED=DE$ and $D^3=E^3=1$ and therefore there is
a unique group homomorphism $\mu:G\rightarrow
\Sym_Z$ such that $\mu(G')=1$, $\mu(d)=D$ and
$\mu(e)=E$. We will check (\ref{IYBS}) for all
$x,y\in Z$. Let $x=m_1d^{i_1}e^{j_1}$ and
$y=m_2d^{i_2}e^{j_2}$, with $m_1,m_2\in G'$,
$i_1,i_2\in\{ 0,1,2\}$ and $j_1,j_2\in\{ 0,1\}$. We
have
\begin{eqnarray*}
x\mu (x)^{-1}(y)&=&
m_1d^{i_1}e^{j_1}\mu(m_1d^{i_1}e^{j_1})\inv(m_2d^{i_2}e^{j_2})\\
&=& m_1m_2
d^{i_1}e^{j_1}E^{2j_1}D^{2i_1}(d^{i_2}e^{j_2})
\hspace{4.4cm} \mbox{(by (\ref{Dk}) and
(\ref{Ek}))}\\ &=& m_1m_2
d^{i_1}e^{j_1}E^{2j_1}(a^{2i_1(i_2+2i_2j_2)+2i_1(2i_1-1)j_2}
b^{i_1i_2} c^{i_1j_2} d^{i_2} e^{j_2})
\hspace{1.4cm} \mbox{(by (\ref{Dk}))} \\ &=& m_1m_2
d^{i_1}e^{j_1} a^{2i_1(i_2+2i_2j_2)+2i_1(2i_1-1)j_2}
b^{i_1i_2+2i_1j_1j_2}
        c^{i_1j_2} d^{i_2+2j_1j_2} e^{j_2} \quad \mbox{(by  (\ref{Ek}))} \\
&=& m_1m_2 a^{2i_1(i_2+2i_2j_2)+2i_1(2i_1-1)j_2+
i_1^2j_2} b^{i_1i_2}
        c^{i_1j_2} d^{i_1}e^{j_1} d^{i_2+2j_1j_2} e^{j_2} \\
&=& m_1m_2 a^{2i_1(i_2+2i_2j_2)+2i_1(2i_1-1)j_2+
i_1^2j_2+2(i_2+2j_1j_2)(i_2+2j_1j_2-1)j_1+i_1(i_2+2j_1j_2)j_1}
\\
 &&\cdot b^{i_1i_2} c^{i_1j_2+(i_2+2j_1j_2) j_1} d^{i_1+i_2+2j_1j_2} e^{j_1+j_2}
 \hspace{4.4cm} \mbox{(by (\ref{diej}))} \\
&=& m_1m_2 a^{2i_1i_2+i_1i_2j_2 + 2i_1^2j_2+i_1j_2+
    2i_2^2j_1+2i_2j_1^2j_2+2j_1^3j_2^2-2i_2j_1-j_1^2j_2+i_1i_2j_1+2i_1j_1^2j_2} \\
&&\cdot b^{i_1i_2} c^{i_1j_2+i_2j_1+2j_1^2j_2}
d^{i_1+i_2+2j_1j_2} e^{j_1+j_2} \\ &=& m_1m_2
a^{2i_1i_2+i_1i_2j_2 + 2i_1^2j_2+i_1j_2+
    2i_2^2j_1+2i_2j_1j_2+2j_1j_2-2i_2j_1-j_1j_2+i_1i_2j_1+2i_1j_1j_2}  \\
&&\cdot b^{i_1i_2} c^{i_1j_2+i_2j_1+2j_1j_2}
d^{i_1+i_2+2j_1j_2} e^{j_1+j_2}
\hspace{1,6cm}\mbox{(since $j_1^2=j_1$ and
$j_2^2=j_2$)}\\
&=& m_1m_2
a^{2i_1i_2+i_1i_2(j_1+j_2)
+2(i_1^2j_2+i_2^2j_1)+i_1j_2+i_2j_1+2(i_1+i_2)j_1j_2+j_1j_2}
\\
 &&\cdot b^{i_1i_2} c^{i_1j_2+i_2j_1+2j_1j_2} d^{i_1+i_2+2j_1j_2} e^{j_1+j_2}\\
\end{eqnarray*}
This expression is invariant by interchanging $m_1$ and $m_2$, $i_1$ and $i_2$, and $j_1$ and
$j_2$. Hence it follows that $x\mu (x)^{-1}(y)=y\mu (y)^{-1}(x)$. By Theorem~\ref{EquivCond}, $G$
is IYB.
 $\qed$}\end{example}

Let $\mu$ be an IYB morphism of $G$.  Then
$\ker(\mu)$ is abelian. If $x \in G$ and $k\in
\ker(\mu)$ then
    \begin{equation}\label{ker}
    \mu(x)(k)=xx\inv\mu(x\inv)\inv(k) = xk\mu(k)(x\inv)=xkx\inv.
    \end{equation}
This implies, by (\ref{2var}), that if $N$ is a
normal subgroup of $G$ contained in $\ker(\mu)$ then
$\mu(x)(Ny)=N\mu(x)(y)$ and therefore (with the
usual bar notation) the map
$\overline{\mu}:\overline{G}=G/N\rightarrow
S_{\overline{G}}$ given by
    $$\overline{\mu}(\overline{x})({\overline{y}}) = \overline{\mu(x)(y)}$$
is well defined and it is easy to show that
$\overline{\mu}$ is an IYB morphism of
$\overline{G}$. We say that $\mu$ is a lifting of
$\overline{\mu}$.

It is somehow natural to try to prove that every
solvable group is IYB with the following induction
strategy: Let $G$ be a non-trivial solvable group.
Take a non-trivial abelian normal subgroup $N$ of
$G$. Assume, by induction, that $\overline{G}=G/N$
has an IYB morphism $\lambda$ and prove that
$\lambda$ admits a lifting to $G$, i.e.
$\lambda=\overline{\mu}$ for some IYB morphism $\mu$
of $G$.  Notice that if this strategy works then
every non-trivial solvable group should have a non
injective IYB morphism. This is the case for all the
examples of IYB groups $G$ that we have computed.
 This leads us to the following natural
question:  Does every IYB group admit a
non-injective IYB morphism? In fact, all the IYB
morphisms which appears implicitly in the results of
Section~\ref{ClassIYB} or in the above examples are
non-injective and therefore they are liftings of IYB
morphisms of proper quotients. These examples may
lead to the impression that every IYB morphism of
$G/N$, for an abelian normal subgroup $N$ of $G$,
can be lifted to an IYB morphism of $G$. However
this is false as the following example shows.

\begin{example}\label{NonLiftable}{\rm
Let $G$ be the group of Example~\ref{C35}. We claim
that the trivial IYB morphism of $G/G'$ does not
lift to an IYB morphism of $G$. Indeed, assume that
$\mu$ is an IYB morphism of $G$ which lifts the
trivial IYB morphism of $G/G'$. Then, clearly,
$\mu(d)$ and $\mu(e)$ commute.  Furthermore, by
(\ref{ker}), $\mu(g)(n)=gng\inv$ for every $g\in G$
and $n\in G'$ and there are $x_j\in G'$,
($j=1,2,3$), with
    \begin{equation}\label{mu35}
    \mu(d)(d)=x_1 d, \quad \mu(e)(d)=x_2 d, \quad \mu(e)(e)=x_3 e.
    \end{equation}
Then, by (\ref{2var}),
    $$dx_2d\inv x_1 d = \mu(d)(x_2d) = \mu(d)\mu(e)(d) = \mu(e)\mu(d)(d)=\mu(e)(x_1d) = ex_1e\inv x_2 d.$$
Therefore $(e,x_1)=(d,x_2)\in \GEN{a}\cap\GEN{b}=\{
1\}$ and hence $x_1,x_2\in Z(G)$. This implies,  by
(\ref{2var}) and (\ref{ker}), that
    $\mu(d^j)(d^k) = x_1^{jk}d^k$ and $\mu(e^j)(d^k) = x_2^{jk}d^k$, for every $j,k$.
Then
$\mu(d)(e)=dd\inv\mu(d\inv)\inv(e)=de\mu(e)\inv(d\inv)=dex_2d^2=x_2c^2e$.
Thus
    $$\mu(e)\mu(d)(e)=\mu(e)(x_2c^2e) = x_2 b^2 c^2 x_3 e$$
and
    $$\mu(d)\mu(e)(e) = \mu(d)(x_3 e) = dx_3d\inv x_2 c^2 e.$$
Therefore $b^2 x_3 = dx_3d\inv$ and so
$b^2=(d,x_3)\in \GEN{a}$, a contradiction.
 $\qed$}\end{example}

\section{Solutions to the Yang-Baxter equation associated to one IYB group}
\label{extsect}

 Let $X$ be a
finite set. We have seen  in
Section~\ref{SecEquivCond} that if $\lambda\colon
X\rightarrow \Sym_X$ is an IYB  map, then the map
$$r\colon X\times X\longrightarrow X\times X$$
defined by
$r(x,y)=(\lambda(x)(y),\lambda(\lambda(x)(y))^{-1}(x))$,
for all $x,y\in X$, is an involutive non degenerate
set theoretical solution of the Yang-Baxter
equation. We also know, by Theorem~\ref{EquivCond},
that $\langle\lambda(X)\rangle$ is an IYB group. One
can ask how many involutive non degenerate set
theoretical solutions of the Yang-Baxter equation
are associated to the same  IYB group. Note that for
any finite set $X$, the map $\lambda\colon
X\rightarrow \Sym_X$ defined by $\lambda (x)=\id_X$
for all $x\in X$ is an IYB map associated to the
trivial group. Note also that if $\lambda_i\colon
X_i\rightarrow \Sym_{X_i}$ is an IYB map associated
to the IYB group $G_i$, for $i=1,\dots, n$, then it
is easy to see that if $X=\bigcup_{i=1}^nX_i$ is a
disjoint union, then the map $\lambda\colon
X\rightarrow \Sym_X$ defined by
$\lambda(x)=\lambda_i(x)$, for all $x\in X_i$, is an
IYB map associated to direct product
$\prod_{i=1}^nG_i$,  considering $\Sym_{X_i}$
naturally included in $\Sym_X$.  Notice that this
also gives another proof for
Corollary~\ref{directproduct}. Hence, since $\langle
1\rangle\times G\cong G$, one can associate with
each IYB group, in an obvious way, infinitely many
involutive non degenerate set theoretical solutions
of the Yang-Baxter equation. Now we give  a
non-obvious construction of an infinite family of
involutive non degenerate set theoretical solutions
of the Yang-Baxter equation associated to a fixed
IYB group.

\begin{lemma}\label{psi}
Let $X$ be a set and let $\lambda\colon X\rightarrow
\Sym_X$ be a map. Let $\psi\colon \Sym_X\rightarrow
\Sym_{X^2}$ be the map defined by
$$\psi(\tau)(x,y)=(\tau(x),\lambda(\tau(x))^{-1}\tau\lambda(x)(y))$$
for all $\tau\in \Sym_X$ and $x,y\in X$. Then $\psi$
is a monomorphism.
\end{lemma}

 \begin{proof}
Let $\tau_1, \tau_2\in \Sym_X$ and $x,y\in X$. We
have:
\begin{eqnarray*}
\psi(\tau_1\tau_2)(x,y)&=&(\tau_1\tau_2(x),\lambda(\tau_1\tau_2(x))^{-1}\tau_1\tau_2\lambda(x)(y))\\
&=&(\tau_1(\tau_2(x)),\lambda(\tau_1(\tau_2(x)))^{-1}\tau_1\lambda(\tau_2(
x)) (\lambda(\tau_2(x))^{-1}\tau_2\lambda(x)(y)))\\
&=&
\psi(\tau_1)(\tau_2(x),\lambda(\tau_2(x))^{-1}\tau_2\lambda(x)(y))\\
&=& \psi(\tau_1)\psi(\tau_2)(x,y).
\end{eqnarray*}
Hence $\psi$ is a homomorphism. It is easy to see
that it is injective.
 \end{proof}

\begin{lemma}\label{product}
Let $\lambda\colon X\rightarrow \Sym_X$ be an IYB
map. Let $\mu\colon X^2\rightarrow \Sym_X$ be the
map defined by $\mu(x,y)=\lambda(x)\lambda(y)$ for
all $x,y\in X$.
 Let $\lambda_2\colon X^2\rightarrow \Sym_{X^2}$ be the map defined by
 $\lambda_2=\psi\mu$, where $\psi$ is as in Lemma~\ref{psi}.
 Then $\lambda_2$ is an IYB map.
 \end{lemma}

 \begin{proof}
 Note that

$$\lambda_2(x,y)\inv(z,t)=(\lambda(y)^{-1}\lambda(x)^{-1}(z),\lambda(\lambda(y)^{-1}\lambda(x)^{-1}(z))^{-1}
 \lambda(y)^{-1}\lambda(x)^{-1}\lambda(z)(t)),$$
 for all $x,y,z,t\in X$.

 We have
 $$\begin{array}{l}
\lambda_2(x,y)\lambda_2(\lambda_2(x,y))\inv(z,t)\\
\hspace{2cm}=
\psi(\lambda(x)\lambda(y)\lambda(\lambda(y)^{-1}\lambda(x)^{-1}(z))\lambda(\lambda(\lambda(y)^{-1}\lambda(x)^{-1}(z))^{-1}
 \lambda(y)^{-1}\lambda(x)^{-1}\lambda(z)(t))).
 \end{array}
 $$

 So, in order to prove the lemma, we need to show that

 \begin{eqnarray}\label{1lemma}
&&\hspace{-2cm}\lambda(x)\lambda(y)\lambda(\lambda(y)^{-1}\lambda(x)^{-1}(z))\lambda(\lambda(\lambda(y)^{-1}\lambda(x)^{-1}(z))^{-1}
 \lambda(y)^{-1}\lambda(x)^{-1}\lambda(z)(t))\notag\\
 &&=
\lambda(z)\lambda(t)\lambda(\lambda(t)^{-1}\lambda(z)^{-1}(x))\lambda(\lambda(\lambda(t)^{-1}\lambda(z)^{-1}(x))^{-1}
 \lambda(t)^{-1}\lambda(z)^{-1}\lambda(x)(y)).
 \end{eqnarray}

 Since $\lambda$ is an IYB map,
 \begin{eqnarray*}
 \lambda(x)\lambda(y)\lambda(\lambda(y)^{-1}\lambda(x)^{-1}(z))&=&
\lambda(x)\lambda(\lambda(x)^{-1}(z))\lambda(\lambda(\lambda(x)^{-1}(z)
)^{-1}(y))\\ &=&
\lambda(z)\lambda(\lambda(z)^{-1}(x))\lambda(\lambda(\lambda(x)^{-1}(z)
)^{-1}(y)).
 \end{eqnarray*}

 Hence (\ref{1lemma}) is equivalent to

 \begin{eqnarray}\label{2lemma}
&&\lambda(\lambda(\lambda(x)^{-1}(z))^{-1}(y))\lambda(\lambda(\lambda(y)^{-1}\lambda(x)^{-1}(z))^{-1}
 \lambda(y)^{-1}\lambda(x)^{-1}\lambda(z)(t))\notag\\
 &&=
 \lambda(\lambda(\lambda(z)^{-1}(x))^{-1}(t))
 \lambda(\lambda(\lambda(t)^{-1}\lambda(z)^{-1}(x))^{-1}
 \lambda(t)^{-1}\lambda(z)^{-1}\lambda(x)(y)).
 \end{eqnarray}

 Since $\lambda$ is an IYB map, to prove (\ref{2lemma}), it is sufficient to
see that

 \begin{eqnarray}\label{3lemma}
\lambda(\lambda(\lambda(x)^{-1}(z))^{-1}(y))(\lambda(\lambda(y)^{-1}\lambda(x)^{-1}(z))^{-1}
 \lambda(y)^{-1}\lambda(x)^{-1}\lambda(z)(t))=
 \lambda(\lambda(z)^{-1}(x))^{-1}(t),
 \end{eqnarray}

 Since $\lambda$ is an IYB map,

 \begin{eqnarray*}
&&\lambda(\lambda(\lambda(x)^{-1}(z))^{-1}(y))\lambda(\lambda(y)^{-1}\lambda(x)^{-1}(z))^{-1}
 \lambda(y)^{-1}\lambda(x)^{-1}\lambda(z)(t)\\
 &&=
 \lambda(\lambda(x)^{-1}(z))^{-1}\lambda(x)^{-1}\lambda(z)(t)\\
 &&=\lambda(\lambda(z)^{-1}(x))^{-1}(t).
 \end{eqnarray*}
 Hence (\ref{3lemma}) is true and the lemma is proved.
 \end{proof}

Note that, with the notation of Lemma~\ref{product},
 if $1\in \lambda(X)$ then $\langle
\lambda(X)\rangle\cong\langle
\lambda_2(X^2)\rangle$. Thus in this way we can
construct infinitely many involutive non degenerate
set theoretical solutions of the Yang-Baxter
equation associated to the IYB group $\langle
\lambda(X)\rangle$.

It seems a difficult problem to describe all the
involutive non degenerate set theoretical solutions
of the Yang-Baxter equation associated to a given
IYB group.

 \vspace{30pt}\small
 \noindent \begin{tabular}{lll}
 F. Ced\'{o} & E. Jespers & \'{A}. del R\'{\i}o \\
 Departament de Matem\`{a}tiques &  Department of Mathematics  & Departamento de Matem\'{a}ticas\\
 Universitat Aut\`{o}noma de Barcelona &  Vrije Universiteit Brussel & Universidad de Murcia \\
 08193 Bellaterra (Barcelona), Spain    &  Pleinlaan 2, 1050 Brussel, Belgium & 30100 Murcia, Spain
\end{tabular}


\begin{thebibliography}{99}
\itemsep=-2pt
\bibitem{BDG} N. Ben David and Y. Ginosar, On groups of central type, non-degenerate and bijective cohomology
 classes, Israel J. Math. to appear, ArXiv: 0704.2516v1 [math.GR].
\bibitem{CJO} F. Ced\'{o}, E. Jespers and J. Okni\'{n}ski,
 The Gelfand-Kirillov dimension of quadratic algebras satisfying the cyclic
 condition, Proc. AMS 134 (2005), 653-663.
\bibitem{drinfeld} V. G. Drinfeld, On unsolved problems in quantum group theory.
Quantum Groups, Lecture Notes Math. 1510,
Springer-Verlag, Berlin, 1992, 1--8.
\bibitem{EG} P. Etingof and S. Gelaki, A method of construction of finite-dimensional triangular semi-
    simple Hopf algebras, Mathematical Research Letters 5 (1998), 551--561.
\bibitem{ESS} P. Etingof, T. Schedler and A. Soloviev,
  Set-theoretical solutions to the quantum
  Yang-Baxter equation, Duke Math. J. 100 (1999), 169-209.
\bibitem{Gat} T. Gateva-Ivanova, A combinatorial approach to the
set-theoretic solutions of the Yang-Baxter equation,
J. Math. Phys. 45 (2004), 3828--3858.
\bibitem{GIVdB} T. Gateva-Ivanova and M. Van den Bergh,
  Semigroups of $I$-type, J. Algebra 206 (1998), 97-112.
\bibitem{JO} E. Jespers and J. Okni\'{n}ski, Monoids and Group of
  $I$-Type, Algebr. Represent. Theory 8 (2005), 709-729.
\bibitem{JObook} E. Jespers and J. Okni\'{n}ski, {\em Noetherian Semigroup Rings}, Springer, Dordrecht 2007.
\bibitem{kassel} C. Kassel, {\em Quantum Groups}, Graduate Text in Mathematics 155, Springer-Verlag, New York, 1995.
\bibitem{passman} D. S. Passman, {\em Permutation Groups}, Benjamin, New York, 1968.
\bibitem{Robinson} D. K. Robinson, {\em A course in the theory of groups}, second edition, Springer-Verlag,
New York, 1996.
\bibitem{Rump} W. Rump, A decomposition theorem for square-free
unitary solutions of the quantum Yang-Baxter
equation, Adv. Math. 193 (2005), 40--55.
\bibitem{Yang} C.N. Yang, Some exact results for the many-body
problem in one dimension with repulsive
delta-function interaction, Phys. Rev. Lett. 19
(1967), 1312--1315.
 \end{thebibliography}
 \end{document}